\documentclass{amsart}
\usepackage{amssymb,latexsym,amsmath,amsfonts}
\input xy
\xyoption{all}
\usepackage[active]{srcltx}
\usepackage[dvips]{color}

\newcommand{\cal}{\mathcal}

\makeatletter
\renewcommand{\subsection}{\@startsection{subsection}{2}{0mm}{-2mm}{-2mm}{\bf\normalsize}}
\def\sbsn{\subsection{\hspace{-3mm}}}

\def\sbsnt#1{\subsection{#1}}
\makeatother

\newtheorem{formula}{}[section]
\newtheorem{definition}[formula]{Definition}
\newtheorem{corollary}[formula]{Corollary}
\newtheorem{remark}[formula]{Remark}
\newtheorem{lemma}[formula]{Lemma}
\newtheorem{theorem}[formula]{Theorem}

\def\thrm{\begin{theorem}}
\def\thrml#1{\begin{theorem}\label{#1}}
\def\ethrm{\end{theorem}}
\def\rmrk{\begin{remark}}
\def\rmrkl#1{\begin{remark}\label{#1}}
\def\ermrk{\end{remark}}
\def\dfntn{\begin{definition}}
\def\dfntnl#1{\begin{definition}\label{#1}}
\def\edfntn{\end{definition}}
\def\nmrt{\begin{enumerate}}
\def\enmrt{\end{enumerate}}
\def\tm#1{\item[{\rm (#1)}]}
\def\qtn{\begin{equation}}
\def\qtnl#1{\begin{equation}\label{#1}}
\def\eqtn{\end{equation}}
\def\lmm{\begin{lemma}}
\def\lmml#1{\begin{lemma}\label{#1}}
\def\elmm{\end{lemma}}
\def\crllr{\begin{corollary}}
\def\crllrl#1{\begin{corollary}\label{#1}}
\def\ecrllr{\end{corollary}}
\def\css{\begin{cases}}
\def\ecss{\end{cases}}

\def\proof{\noindent{\bf Proof}.\ }

\def\G{{\cal G}}

\def\X{{\cal X}}
\def\Y{{\cal Y}}

\def\J{{\cal J}}

\def\ZZ{{\mathbb Z}}

\def\SD{{\scriptscriptstyle\Delta}}
\def\SG{{\scriptscriptstyle\Gamma}}

\DeclareMathOperator{\aut}{Aut}

\DeclareMathOperator{\Fib}{Fib}

\DeclareMathOperator{\id}{id}

\DeclareMathOperator{\iso}{Iso}
\DeclareMathOperator{\mat}{Mat}
\DeclareMathOperator{\orb}{Orb}

\DeclareMathOperator{\sym}{Sym}

\def\bull{\hfill\vrule height .9ex width .8ex depth -.1ex\medskip}

\def\ddp{{\scriptscriptstyle\bot}}
\def\lg{\langle}

\def\ov{\overline}
\def\rg{\rangle}
\def\wt{\widetilde}

\def\VRT#1{*=<5mm>[o][F-]{#1}}

\def\grphp#1{$\xymatrix@R=10pt@C=10pt@M=0pt@L=2pt{#1}$}
\def\mult{}

\newcommand{\comment}[1]{}

\begin{document}
\title{On quasi-thin association schemes}
\author{Mikhail Muzychuk}
\address{Netanya Academic College, Netanya, Israel}
\email{muzy@netanya.ac.il}
\author{Ilya Ponomarenko}
\address{Steklov Institute of Mathematics at St. Petersburg, Russia}
\email{inp@pdmi.ras.ru}
\date{}

\begin{abstract}
An association scheme is called quasi-thin if the valency of each its basic relation is one
or two. A quasi-thin scheme is Kleinian if the thin residue of it forms a Klein group with
respect to the relation product. It is proved that any Kleinian scheme arises from near-pencil
on~$3$ points, or affine or projective plane of order~$2$. The main result is that any
non-Kleinian quasi-thin scheme a) is the two-orbit scheme of a suitable permutation group,
and b) is characterized up to isomorphism by its intersection number array. An infinite
family of Kleinian quasi-thin schemes for which neither a) nor b) holds is also constructed.
\end{abstract}

\maketitle

\section{Introduction}
Given a permutation group $G\le\sym(\Omega)$ one can define a {\it schurian} coherent configuration
$(\Omega,S)$ where $S$ is the set of $G$-orbits with respect to the componentwise action
of $G$ on the set $\Omega\times\Omega$ (as for a background of association schemes
and coherent configurations see Section~\ref{220309a}). However, not
all coherent configuration can be obtained in this way. This leads naturally to so called
{\it schurity problem}: find an internal characterization of schurian coherent configurations
in a given class. Sometimes a solution of this problem is obtained by proving that any
coherent configuration from the class is {\it separable}, i.e. is characterized
up to isomorphism by the intersection number array. The {\it separability problem} consists
in finding an internal characterization of separable coherent configurations in a given
class. A comprehensive survey of the schurity and separability problems can be found in~\cite{EP09}.\medskip

In this paper we deal with the schurity and separability problems in the class of
quasi-thin association schemes: an association scheme is called {\it quasi-thin} if
the valency of each its basic relation is one or two.\footnote{In \cite{EP00} it was proved
that each {\it thin} or {\it regular} association scheme, i.e. such that the valences
are ones, is schurian and separable.}
Every finite group $G$ of even order with a chosen involution $a$ gives rise to schurian
quasi-thin scheme corresponding the action of $G$ on cosets modulo $\lg a\rg$.
Despite of the fact that quasi-thin schemes were introduced  explicitly only in 2002 (\cite{HM02b}),
the first result about  quasi-thin scheme goes back to~\cite{W76} where it
was proved that any primitive quasi-thin scheme is schurian (and, in fact separable). Only
quarter of century later this result was generalized to some special classes of quasi-thin
schemes \cite{H01,HM02a,MZ}. However, nothing was known on their separability. On the other
hand, there are non-schurian and non-separable quasi-thin schemes: in the Hanaki-Miyamoto list~\cite{HM}
one can find $1$, $1$ and $26$ non-schurian quasi-thin schemes on $16$, $28$ and $32$
points respectively, and the schemes on $16$ and $28$ points are non-separable.
In all these examples the scheme in question has a very "special structure" explained
below.\medskip

Let $\X=(\Omega,S)$ be an association scheme and $S_1$ the set of all valency one
relations in~$S$. We say that $\X$ is a {\it Kleinian} scheme if the thin residue of it
is contained in $S_1$ and forms a Klein group with respect to the relation product.
In this case the degree $n=|\Omega|$ of the scheme~$\X$ is divided by~$4$. Moreover,
one can prove that~$\X$ is the algebraic fusion of a coherent configuration with $n/4$
regular homogeneous components of degree~$4$ by means of a group of algebraic
automorphisms acting regularly on the set of fibers (see Subsection~\ref{030910a}).
The number $n/n_1$ where $n_1=|S_1|$ will be called the {\it index} of the scheme~$\X$.
The following statement is the main result of the paper.

\thrml{210210a}
Any non-schurian or non-separable quasi-thin scheme is a Kleinian scheme of index~$4$ or~$7$.
Moreover, given $i\in\{4,7\}$ there exist infinitely many both non-schurian and non-separable
Kleinian schemes of index~$i$.
\ethrm

The proof is given in Subsection~\ref{030910a} and is divided into two parts. In the
first of them we prove that all non-Kleinian quasi-thin schemes are schurian and separable.
For this purpose we introduce the notion of orthogonal in a quasi-thin scheme and show
that any quasi-thin scheme with at most one orthogonal is schurian and separable
(Section~\ref{090209g}). To deal with the remaining non-Kleinian schemes we study the
one point extension of a quasi-thin scheme (Section~\ref{030910b}) and give
a sufficient condition for such a scheme to be schurian and separable in terms of the
existence of the one point extension of an algebraic isomorphism (Theorem~\ref{250210r}).
The key point in the first part of the proof is Theorem~\ref{280109e} showing that
in our case such extension does always exist.\medskip

The second part of the proof deals with a Kleinian quasi-thin scheme.
Every such a scheme is an algebraic fusion of a coherent configuration each homogeneous
component of which is the scheme of a Klein group. These configurations are called
Kleinian and studied in Section~\ref{290810d}. We show that any such a configuration
is closely related to a partial linear space, and classify all possible spaces
in Corollary~\ref{010710a}. For Kleinian configuration arising from quasi-thin schemes
this classification reduced to three cases: near-pencil on three points and affine or
projective plane of order~$2$ (Corollary~\ref{020710u}). The Kleinian schemes of the
near-pencil type are schurian and separable whereas in the other two cases we
construct infinitely many non-schurian and non-separable quasi-thin schemes.

\crllrl{260810a}
Any non-Kleinian quasi-thin scheme is schurian and separable.\bull
\ecrllr

It would be too naive to expect that any commutative quasi-thin scheme is always schurian
and separable because such a scheme can be Kleinian. Indeed, let $A_1$ be the direct
product of two cyclic groups of order~$4$ and $f_1$ the involutive automorphism of $A_1$
taking $a$ to $a^{-1}$; let $A_2$ be the direct product of two Klein groups and $f_2$ the
involutive automorphism of $A_2$ which interchanges the coordinates. Denote by $\X_i$
the scheme of the permutation group on~$A_i$ generated by the regular representation of
$A_i$ and the automorphism $f_i$, $i=1,2$. Then $\X_1$ and $\X_2$ are commutative schurian
quasi-thin  schemes of degree~$16$ and rank~$10$. Moreover, a direct computation shows that
they are (a) Kleinian, (b) non-isomorphic and (c) algebraically isomorphic. In particular,
none of them is separable. In contrast to this example we prove the following theorem.

\thrml{180510a}
A commutative quasi-thin scheme is schurian.
\ethrm

The proof of this theorem is reduced by Theorem~\ref{210210a} to the case of commutative
Kleinian quasi-thin scheme. The schurity of such a scheme is proved by a direct
computation in Subsection~\ref{030910c}.\medskip

All undefined terms and notation concerning permutation groups can be found in~\cite{DM}.
To make the paper self-contained we give a background on theory of coherent configurations
and on schurity and separability problems in Sections~\ref{220309a} and~\ref{290810c}.\medskip

{\bf Notation.}
Throughout the paper $\Omega$ denotes a finite set. The diagonal of the Cartesian square
$\Omega\times\Omega$ is denoted by~$1_\Omega$; for any $\alpha\in\Omega$ we set
$1_\alpha=1_{\{\alpha\}}$. For a relation $r\subset\Omega\times\Omega$ we set
$r^*=\{(\beta,\alpha):\ (\alpha,\beta)\in r\}$ and
$\alpha r=\{\beta\in\Omega:\ (\alpha,\beta)\in r\}$ for all $\alpha\in\Omega$.
For $\Gamma,\Delta\subset\Omega$ we set $u_{\SG,\SD}=u\cap(\Gamma\times\Delta)$
and $u_\SG=u_{\SG,\SG}$.
For $s\subset \Omega\times\Omega$ we set
$r\cdot s=\{(\alpha,\gamma):\ (\alpha,\beta)\in r,\ (\beta,\gamma)\in s$
for some $\beta\in\Omega\}$. If $S$ and $T$ are sets of relations, we set
$S\cdot T=\{s\cdot t:\ s\in S,\, t\in T\}$. The set of all unions of the elements of $S$
is denote by $S^\cup$.

\section{Association schemes and coherent configurations}\label{220309a}
This section accumulates the basic definitions and facts about coherent configurations
and association schemes which are needed for understanding the paper (see also~\cite{EP09,Zi1}).

\sbsnt{Definitions.}
A pair $\X=(\Omega,S)$ where $\Omega$ is a finite set and $S$ a partition of
$\Omega\times\Omega$, is called a {\it coherent configuration} on $\Omega$ if
$1_\Omega\in S^\cup$, $S^*=S$ and given $u,v,w\in S$, the number
$$
c_{uv}^w=|\alpha u\cap \gamma v^*|
$$
does not depend on the choice of $(\alpha,\gamma)\in w$. The elements of $\Omega$, $S$,
$S^\cup$ and the numbers (S3) are called the {\it points}, the {\it basic relations}, the
{\it relations} and the {\it intersection numbers} of~$\X$, respectively. The numbers
$|\Omega|$ and $|S|$ are called the {\it degree} and {\it rank} of it. The coherent
configuration $\X$ is {\it commutative} if $c_{uv}^w=c_{vu}^w$ for all $u,v,w\in S$.
The unique basic relation containing a pair $(\alpha,\beta)\in\Omega\times\Omega$ is denoted
by $r(\alpha,\beta)$.\medskip

For the intersection numbers we have the following well-known
identities (see~\cite{Hig75}).
$$
c_{u^*v^*}^{w^*}=c_{vu}^w\quad\text{and}\quad |w|c_{uv}^{w^*}=|u|c_{vw}^{u^*}=|v|c_{wu}^{v^*},\qquad u,v,w\in S
$$
If the configuration is homogeneous (a scheme), then these equalities my be rewritten as follows
(see~\cite{Zi1}):
\qtnl{150410a}
c_{u^*v^*}^{w^*}=c_{vu}^w\quad\text{and}\quad n_wc_{uv}^{w^*}=n_uc_{vw}^{u^*}=n_vc_{wu}^{v^*},\qquad u,v,w\in S
\eqtn
 The set of basic relations contained in $u\cdot v$ with $u,v\in S^\cup$ is
denoted by~$uv$. Sometimes it is useful to treat $uv$ as a multiset in which an element
$w\in S$ appears with the multiplicity $c_{uv}^w$. This multiset will be written
as the element of the free module $\ZZ S$ equipped by the involution~$*$ and the natural
scalar product defined by
$$
(\sum_{s\in S}x_ss)^*=\sum_{s\in S}x_ss^*\quad\text{and}\quad
\lg\sum_{s\in S}x_ss,\sum_{s\in S}y_ss\rg=\frac{1}{|\Omega|}\sum_{s\in S}x_sy_s|s|.
$$
Notice that this a scalar product is associative, that is
\qtnl{associative}
\lg xy,z\rg=\lg y,x^*z\rg,\qquad x,y,z\in\ZZ S.
\eqtn
For a homogeneous configuration the scalar product reads as follows:
$$
\lg\sum_{s\in S}x_ss,\sum_{s\in S}y_ss\rg=\sum_{s\in S}x_sy_s n_s.
$$

Each time we use notation $uv$ it will be clear is it a set or multiset.

\sbsnt{Fibers and homogeneity.}
Any set $\Delta\subset\Omega$ for which $1_\Delta\in S$, is called the {\it fiber} of the coherent
configuration~$\X$; the set of all of them is denoted by $\Fib(\X)$. Clearly, the set
of points is the disjoint union of fibers. One can also see that if $\Delta$ is a union of
fibers and $S_\Delta$ is the set of all nonempty relations $u_\SD$ with
$u\in S$, then $(\Delta,S_\Delta)$ is a coherent configuration, called the {\it restriction}
of~$\X$ to~$\Delta$. Besides, for any basic
relation $u\in S$ there exist uniquely determined fibers $\Delta,\Gamma$ such that
$u\subset\Delta\times\Gamma$. Set
$$
n_u=c_{uu^*}^v
$$
where $v=1_\Delta$. Then the number $|\delta u|=n_u$ does not depend on $\delta\in\Delta$.
When $n_{u^{}}=n_{u^*}=1$, the relation $u$ is called {\it thin}.\medskip

The coherent configuration $\X$ is called {\it homogeneous} or a {\it scheme} if
$1_\Omega\in S$, or equivalently if $\Fib(\X)=\{\Omega\}$. In this case $n_u=n_{u^*}$
for all $u\in S$; the number $n_u$ is called the {\it valency} of~$u$. The set of all
basic relations of valency~$m$ is denoted by~$S_m$. The following result proved in~\cite{MP} will
be used in Section~\ref{090209g}.

\lmml{141208j}
Let $(\Omega,S)$ be a scheme and $u,v\in S$. Then $c_{u^*v}^w\le 1$ for all
$w\in S$ if and only if $uu^*\cap vv^*=\{1_\Omega\}$.\bull
\elmm

\sbsnt{Closed sets.} Let $\X=(\Omega,S)$ be a scheme. A set $T\subset S$ is
called {\it closed}, notation $T\leq S$, if $TT^*\subset T$. It is easily seen that the set $S_1$, called the {\it thin radical}of~$\X$ in \cite{Zi1}, is closed
(and forms a group with respect to the relational product).
The intersection of all closed sets containing  the set $\bigcup_{u\in S}uu^*$ is called the {\it thin residue} of~$\X$.
The union of all relations from a closed set~$T$ is an equivalence relation on~$\Omega$ with
classes $\alpha T$, $\alpha\in\Omega$. The set of all these classes is denoted by $\Omega/T$.
One can prove that the pairs
$$
\X_\Delta=(\Delta,S_\Delta)\quad\text{and}\quad\X_{\Omega/T}=(\Omega/T,S_{\Omega/T})
$$
are schemes where $\Delta\in\Omega/T$ and $S_\Delta$ is as above, and $S_{\Omega/T}$ consists
of all relations of the form $\{(\Delta,\Gamma)\in\Omega/T\times\Omega/T:\ u_{\SD,\SG}\ne\emptyset\}$
with $u\in S$. The schemes $\X_\Delta$ and $\X_{\Omega/T}$ are called the {\it restriction}
of $\X$ to~$\Delta$, and the {\it quotient} of~$\X$ modulo~$T$.

\sbsnt{Extensions.} There is a natural partial order\, $\le$\, on the set of all coherent
configurations on the set~$\Omega$. Namely, given two coherent configurations
$\X=(\Omega,S)$ and $\X'=(\Omega,S')$ we set
$$
\X\le\X'\ \Leftrightarrow\ S^\cup\subset (S')^\cup.
$$
In this case $\X'$ is called an {\it extension} or {\it fission} of~$\X$. The minimal and
maximal elements with respect to that order are respectively the coherent configurations of
rank~$2$ and of rank $n^2$ where $n=|\Omega|$. The first of them is called {\it trivial}; its
basic relations are $1_\Omega$ and $\Omega\times\Omega\setminus\{1_\Omega\}$. The second
one is called {\it complete}; in this case $S^\cup$ consists of all binary relations
on~$\Omega$.\medskip

Let $\X=(\Omega,S)$ be a coherent configuration and $\alpha\in\Omega$. Denote by $S_\alpha$
the set of basic relations of the smallest coherent configuration on $\Omega$ such that
$$
1_{\alpha}\in S_\alpha\quad\text{and}\quad S\subset S_\alpha^\cup.
$$
The coherent configuration $\X_\alpha=(\Omega,S_\alpha)$ is called the {\it $\alpha$-extension}
(or a {\it one point extension}) of the coherent configuration~$\X$. It is easily seen that
given $u,v,w\in S$ the set $\alpha u$ and the relation $w_{\scriptstyle\alpha u,\alpha v}$
are unions of some fibers and some basic relations of the coherent configuration~$\X_\alpha$,
respectivelt.\medskip

Let $\X=(\Omega,S)$ be a scheme and $T\subset S$ a closed set containing the thin residue
of~$\X$. Denote by $S_{(T)}$
the set of all relations $u_{\SD,\SG}$ where $u\in S$ and $\Delta,\Gamma\in\Omega/T$.
Then from~\cite[Theorem~2.1]{EP09a} (see also \cite{28}) it follows that the pair
$\X_{(T)}=(\Omega,S_{(T)})$
is a coherent configuration; it is called the {\it thin residue extension} of the scheme~$\X$.

\sbsnt{$1$-regularity.}\label{250210c} Let $\X=(\Omega,S)$ be a coherent configuration.
A point $\alpha\in\Omega$ is called {\it regular} (in~$\X$), if
\qtnl{040409b}
|\alpha u|\le 1,\qquad u\in S.
\eqtn
Suppose that the set $\Delta$ of all regular points is nonempty. Then the coherent configuration~$\X$
is called {\it $1$-regular}. In this case all basic relations of the coherent configuration
$\X_\Delta$ are thin. A $1$-regular scheme is called {\it regular}; regular schemes are
exactly thin schemes in the sense of~\cite{Zi1}. One can prove that if $\X$ is a scheme
and $T$ is the thin residue of~$\X$, then the scheme $\X_{\Omega/T}$ is regular.

\sbsnt{Direct sum and tensor product.}
Let $\X=(\Omega,S)$ and $\X'=(\Omega',S')$ be two coherent configurations. Denote by
$\Omega\sqcup\Omega'$ the disjoint union of~$\Omega$ and~$\Omega'$, and by
$S\boxplus S'$ the union of the set $S\cup S'$ and the set of all relations
$\Delta\times\Delta'$ and $\Delta'\times\Delta$ with $\Delta\in\Fib(\X)$ and $\Delta'\in\Fib(\X')$. Then the pair
$$
\X\boxplus\X'=(\Omega\sqcup\Omega',S\boxplus S')
$$
is a coherent configuration called the {\it direct sum} of~$\X$ and~$\X'$.
One can see that $\X\boxplus\X'$ is the smallest coherent configuration
(on~$\Omega\sqcup\Omega'$) the restriction of which to $\Omega$ and $\Omega'$
are respectively~$\X$ and~$\X'$. It should be noted that the direct sum of
any two coherent configurations is non-homogeneous.\medskip

Set $S\otimes S'=\{u\otimes u':\ u\in S,\ u\in S'\}$ where $u\otimes u'$ is
the relation on $\Omega\times\Omega'$ consisting of all pairs $((\alpha,\alpha'),(\beta,\beta'))$
with $(\alpha,\beta)\in u$ and $(\alpha',\beta')\in u'$. Then the pair
$$
\X\otimes\X'=(\Omega\times\Omega',S\otimes S')
$$
is a coherent configuration called the {\it tensor product} of~$\X$ and~$\X'$.
It should be noted that it is homogeneous if only if so are the factors.

\section{Schurian and separable coherent configurations}\label{290810c}

\sbsnt{Isomorphisms and schurity.}
Two coherent configurations are called {\it isomorphic} if there exists a bijection between
their point sets preserving the basic relations. Any such bijection is called the
{\it isomorphism} of these coherent configurations. The group of all isomorphisms of a
coherent configuration $\X=(\Omega,S)$ contains a normal subgroup
$$
\aut(\X)=\{f\in\sym(\Omega):\ u^f=u,\ u\in S\}
$$
called the {\it automorphism group} of~$\X$. It is easily seen that given $\alpha\in\Omega$
we have $\aut(\X)_\alpha=\aut(\X_\alpha)$ where $\X_\alpha=(\Omega,S_\alpha)$.\medskip

Conversely, let $G\le\sym(\Omega)$ be a permutation group and $S$ the set of orbits of the
componentwise action of $G$ on~$\Omega\times\Omega$. Then $\X$ is a coherent configuration
and we call it the {\it coherent configuration of~$G$}. This coherent configuration is
homogeneous if and only if the group is transitive; in this case we say that $\X$ is the
{\it scheme of~$G$}. A coherent configuration on $\Omega$ is called {\it schurian} if it is
the coherent configuration of some permutation group on~$\Omega$. It is easily seen that
a coherent configuration~$\X$ is schurian if and only if it is the coherent configuration
of the group~$\aut(\X)$.

\sbsnt{Algebraic isomorphisms and separability.}
Two coherent configurations $\X=(\Omega,S)$ and $\X'=(\Omega',S')$ are called
{\it algebraically isomorphic} if
\qtnl{f041103p1}
c_{uv}^w=c_{u'v'}^{w'},\qquad u,v,w\in S,
\eqtn
for some bijection $\varphi:S\to S',\ u\mapsto u'$ called the {\it algebraic isomorphism}
from~$\X$ to~$\X'$. Each isomorphism $f$ from~$\X$ to~$\X'$
induces in a natural way an algebraic isomorphism between these schemes denoted by $\varphi_f$.
The set of all isomorphisms inducing the algebraic isomorphism~$\varphi$ is denoted by
$\iso(\X,\X',\varphi)$. In particular,
$$
\iso(\X,\X,\id_S)=\aut(\X)
$$
where $\id_S$ is the identical mapping on $S$. A coherent configurations $\X$ is
called {\it separable} if for any algebraic isomorphism~$\varphi:\X\to\X'$
the set $\iso(\X,\X',\varphi)$ is a non-empty one. Given points $\alpha\in\Omega$
and $\alpha'\in\Omega'$ an algebraic isomorphism $\varphi':\X^{}_{\alpha^{}}\to\X'_{\alpha'}$
is called the {\it $(\alpha,\alpha')$-extension} or the {\it one point extension}
of~$\varphi$ if
\qtnl{290309z}
\varphi'(1_{\alpha^{}})=1_{\alpha'},\qquad \varphi'(u)\subset\varphi(\wt u),\ u\in S_\alpha,
\eqtn
where $\wt u$ is the unique basic relation of $\X$ that contains $u$. Clearly,
$\varphi'$ is uniquely determined by $\varphi$.

\sbsnt{Examples.}
One can see that a coherent configuration is $1$-regular if and only if it is a coherent
configuration of a permutation group having a faithful regular orbit. The proof of this
statement as well as the next one can be found in~\cite{EP00}.

\thrml{210109e}
Any $1$-regular coherent configuration is schurian and separable.\bull
\ethrm

One can define the class of all coherent configurations that can be constructed from
$1$-regular coherent configurations by means of direct sums and tensor products.
By Theorem~\ref{210109e} and the following statement proved in \cite[Theorems~1.17,1.20]{E05}
any coherent configuration from this class is schurian and separable.

\thrml{240810a}
Let $\X_1$, $\X_2$ be coherent configurations and let $\X$ be $\X_1\boxplus\X_2$ or
$\X_1\otimes\X_2$. Then $\X$ is schurian (resp. separable) if and only if both $\X_1$
and $\X_2$ are schurian (resp. separable).\bull
\ethrm

\sbsnt{Thin residue extension.} In this section we study the schurity and separability
of thin residue extension of an arbitrary scheme.

\thrml{150610b}
Any scheme with the separable thin residue extension is separable, and is schurian if and
only if so is the extension.
\ethrm
\proof Let $\X=(\Omega,S)$ be a scheme and $T$ its thin residue. Then from~\cite[Theorem~2.1]{EP09a}
it follows that the following statements hold:\medskip

\noindent (i) given $f\in\aut(\X_{\Omega/T})$ the mapping
$\psi_f:u_{\SD^{},\SG^{}}\mapsto u_{\SD^f,\SG^f}$ is
an algebraic isomorphism of the coherent configuration $\X_{(T)}$,\footnote{This
statement was also proved in \cite{28}.}\vspace{2mm}

\noindent (ii) 
given an algebraic isomorphism $\varphi:\X\to\X'$ and
$f\in\iso(\X^{}_{\Omega^{}/T^{}},\X'_{\Omega'/T'},\varphi^{}_{\Omega/T})$ where $T'=T^\varphi$,
there exists an algebraic isomorphism $\varphi_T:\X^{}_{(T^{})}\to\X'_{(T')}$
extending $\varphi$ and such that $(\varphi_T)^{}_{\Omega/T}$ is induced by~$f$.\medskip

Suppose that the thin residue extension $\X_0=\X_{(T)}$ of the scheme~$\X$ is separable. Then
the first statement immediately follows from statement~(ii). To prove the second one suppose
first that the coherent configuration $\X_0$ is schurian. Take a relation
$u\in S$ and pairs $(\alpha_1,\beta_1)$,
$(\alpha_2,\beta_2)\in u$. Set
$$
\Delta_i=\alpha_i T\quad\text{and}\quad \Gamma_i=\beta_i T,\qquad i=1,2.
$$
Then the pairs $(\Delta_1,\Gamma_1)$ and $(\Delta_2,\Gamma_2)$ belong to the same
relation of the quotient scheme $\X_{\Omega/T}$. Since this scheme is regular (and hence schurian),
one can find an automorphism $f\in\aut(\X_{\Omega/T})$ taking $(\Delta_1,\Gamma_1)$ to
$(\Delta_2,\Gamma_2)$. By statement~(i) it induces the algebraic
isomorphism $\psi_f$ of the coherent configuration~$\X_0$ such that
\qtnl{150610c}
(\Delta_1)^{\psi_f}=\Delta_2\quad\text{and}\quad(\Gamma_1)^{\psi_f}=\Gamma_2.
\eqtn
Since this coherent configuration is separable, the algebraic isomorphism
$\psi_f$ is induced by an isomorphism~$g$ of~$\X_0$ to itself. From the definition of $\psi_f$
it follows that $g\in\aut(\X)$. Moreover, due to \eqref{150610c} we also have
$$
(\Delta_1)^g=\Delta_2\quad\text{and}\quad(\Gamma_1)^g=\Gamma_2.
$$
Thus without loss of generality we can assume that $\Delta_1=\Delta_2$ and
$\Gamma_1=\Gamma_2$. Denote these sets by $\Delta$ and $\Gamma$. Then the pairs
$(\alpha_1,\beta_1)$ and $(\alpha_2,\beta_2)$ belong to the relation $u_{\Delta,\Gamma}$
and we are done by the schurity of~$\X_0$.\medskip

To complete the proof suppose that the scheme $\X$ is schurian. Take a basic relation $u_0$
of the coherent configuration~$\X_0$ and pairs $(\alpha,\beta)$, $(\alpha',\beta')\in u_0$.
Then $u_0$ is contained in a certain relation $u\in S$, and
$$
(\alpha,\beta),\ (\alpha',\beta')\in u\quad\text{and}\quad
\alpha T=\alpha' T,\ \beta T=\beta' T.
$$
Denote the latter two sets (which are elements of $\Omega/T$) by $\Delta$ and $\Gamma$.
By the schurity of $\X$ one can find $f\in\aut(\X)$ such that
$(\alpha^f,\beta^f)=(\alpha',\beta')$. Clearly, $\Delta^f=\Delta$ and
$\Gamma^f=\Gamma$. On the other hand, since the scheme $\X_{\Omega/T}$ is regular, we have
$$
\aut(\X_0)=\aut(\X)_{\{\Delta\}}\cap \aut(\X)_{\{\Gamma\}}.
$$
Thus $f\in\aut(\X_0)$, and the coherent configuration $\X_0$ is schurian.\bull

The assumption on the separability of the thin residue extension in
Theorem~\ref{150610b} is essential. Indeed, let $\X$ be the quasi-thin scheme of degree~$16$
that has number \#173 in~\cite{HM}. Then the group $\aut(\X)$ has two orbits, and hence
the scheme $\X$ is non-schurian. On the other hand, its thin residue extension
is schurian but non-separable.

\section{Klein configurations}\label{290810d}

\sbsnt{Definition and structure.}\label{010710s}
Throughout this section we fix a Klein group~$G$. A coherent configuration $\X=(\Omega,S)$
is a {\it Klein configuration} if any its homogeneous component is the scheme of a regular
permutation group isomorphic to~$G$. In this case we fix a semiregular action of~$G$ on $\Omega$
such that $\orb(G,\Omega)=\Fib(\X)$ and the homogeneous component $\X_\Delta$ corresponding to
a fiber $\Delta\in\Fib(\X)$ is the scheme of the group $G^\Delta$. This semiregular action
of~$G$ is completely determined by the group isomorphisms
\qtnl{240610b}
G\to S^{}_{\Delta,\Delta},\ g\mapsto g^{}_\Delta
\eqtn
where given $\Gamma,\Delta\in\Fib(\X)$ we set $S_{\Gamma,\Delta}=\{s\in S:\ s\subset\Gamma\times\Delta\}$.
It should be noted that these isomorphisms, and hence the semiregular action of the
group~$G$, can be chosen not a unique way.\medskip

Let $\X=(\Omega,S)$ be a Klein configuration and $\Fib(\X)=\{\Omega_i\}_{i\in I}$ where
$I$ is a nonempty finite set. Then due to~\cite[Lemma~5.1]{EP99} given indices $i,j\in I$
the groups
\qtnl{020710a}
L_{ij}=\{g\in G:\ g_i\cdot s=s\},\quad R_{ij}=\{g\in G:\ s\cdot g_j=s\}
\eqtn
do not depend on the relation $s\in S_{ij}$ where $g_i=g^{}_{\Omega_i}$,
$g_j=g^{}_{\Omega_j}$ and $S_{ij}=S_{\Omega_i,\Omega_j}$. Moreover, from the same result
it follows that $L_{ij}=R_{ji}$, $R_{ij}=L_{ji}$ and exactly one of the following three
statements hold:
\nmrt
\tm{K1} $L_{ij}=R_{ij}=G$ and $S_{ij}=\{\Omega_i\times\Omega_j\}$,
\tm{K2} $L_{ij}=R_{ij}=\{1\}$ and $S_{ij}=\orb(G,\Omega_i\times\Omega_j)$,
\tm{K3} $|L_{ij}|=|R_{ij}|=2$ and $S_{ij}=
\{\Omega_{i,1^{}}\times\Omega_{j,1^f}\,\cup\, \Omega_{i,2^{}}\times\Omega_{j,2^f}\}_{f\in\sym(2)}$
\enmrt
where $\{\Omega_{i,1},\Omega_{i,2}\}=\orb(L_{ij},\Omega_i)$ and
$\{\Omega_{j,1},\Omega_{j,2}\}=\orb(R_{ij},\Omega_j)$. In what follows the
array $R=R(\X,G)=(R_{ij})$ is treated as a matrix whose rows and columns are indexed by the
elements of the set~$I$. Clearly, given $i,j,k\in I$ we have
\qtnl{300610c}
R_{ji}=R_{ki}\quad\&\quad|R_{jk}|=|R_{ji}|=2\qquad\Rightarrow\qquad R_{ij}=R_{kj}\quad\&\quad R_{ik}=R_{jk}.
\eqtn
(Indeed, we always have $S_{ji}\cdot S_{ik}\subset S_{jk}^\cup$, and the left-hand side
conditions imply that in fact $S_{ji}\cdot S_{ik}=S_{jk}$.) Moreover, the relation~$\sim$
consisting of all pairs $(i,j)\in I\times I$ such that $R_{ij}=\{1\}$, is an equivalence relation, and
\qtnl{300610z}
j\sim k\ \Rightarrow\ R_{ji}=R_{ki},\qquad i\in I.
\eqtn
It should be noted that the matrix $R$ depends on the choice of isomorphisms~\eqref{240610b}.
On the other hand, the conditions (K1), (K2) and (K3) imply the following statement.

\lmml{240610c}
Let $\X$ and $\X'$ be Klein configurations on the same set. Suppose that $\Fib(\X)=\Fib(\X')$
and $R(\X,G)=R(\X',G)$. Then $\X=\X'$.\bull
\elmm

\sbsnt{Reduced configurations.} Given a set $J\subset I$ denote by $\X_J$ the restriction of
the Klein configuration~$\X$ to the union $\Omega_J$ of all fibers $\Omega_i$, $i\in J$. Then
$\X_J$ is also a Klein configuration and $R_J=R(\X_J,G)$ is a submatrix of~$R$ the rows and
columns of which are the elements of~$J$. Denote by $\J(\X,G)$ the set of all transversals of
the equivalence relation $\sim$.

\lmml{290610k}
Given $J\in\J(\X,G)$ the Klein configurations $\X$ and $\X_J$ are schurian (or separable)
simultaneously.
\elmm
\proof By the lemma hypotheses for any fiber $\Omega_i$ with $i\in I\setminus J$
any relation of the set $S_{ij}$ with $j\sim i$ is thin. Thus the required statement
immediately follows from statement~(2) of \cite[Lemma~9.4]{EP02}.\bull

The Klein configuration $\X_J$ from Lemma~\ref{290610k} is {\it reduced}: by the
definition this means that the equivalence relation $\sim$ is trivial, or equivalently
$|R_{ij}|\ge 2$ for all distinct $i,j\in I$. For any reduced configuration we define the
incidence structure $\G=(I,L)$ with the point set $I$ and the line set $L$ consisting
of all sets
$$
L_i(H)=\{i\}\,\cup\,\{j\in I:\ R_{ji}=H\}
$$
where $H=R_{ji}$ for some $j\in I\setminus\{i\}$. (Here $H$ is always a subgroup of
$G$ of order~$2$.) In the following statement we show that under rather a technical
assumption the geometry $\G$ is a partial linear space, i.e. an incidence
structure such that the lines have size at least~$2$ and two distinct points are incident
to at most one line.

\lmml{290610c}
Let $\X$ be a reduced Klein configuration such that for any $i\in I$ there exists $j\in I\setminus\{i\}$
with $R_{ji}\ne G$. Then $\G$ is a partial linear space in which any point is incident to at
most three lines.
\elmm
\proof By the hypothesis any line $L_i(H)\in L$ contains at least two points: $i$
and $j\in I\setminus\{i\}$ for which $R_{ji}\ne G$. Next, from~\eqref{300610c} it follows
that given an element $i\in I$ and a group $H\le G$ of order~$2$ we have
\qtnl{290610y}
j\in L_i(H)\ \Rightarrow\ L_j(K)=L_i(H)
\eqtn
where $K=R_{ij}$. Next, suppose that distinct points $i$ and $j$ are incident to two lines
$L_k(H)$ and $L_{k'}(H')$ where $k,k'\in I$ and $H,H'\le G$ are of order~$2$. Then due
to~\eqref{300610c} we have $R_{ki}=R_{ji}=R_{k'i}$. Denote this group by $K$. Clearly,
$|K|=2$. Therefore due to (\ref{290610y}) we conclude that
$L_k(H)=L_i(K)=L_{k'}(H')$. Thus any two distinct points are incident to at most one line.
Since the group $G$ has exactly three subgroups of order~$2$, any point is incident to at
most three lines.\bull

We recall that a linear space is partial linear space in which any two distinct points
are incident to exactly one line.

\crllrl{010710a}
In the condition of Lemma~\ref{290610c} suppose that $R_{ij}\ne G$ for all $i,j\in I$. Then
either $|L|=1$, or $|I|\le 7$ and $\G$ is a projective or affine plane of order~$2$, or
$\G$ is one of the four linear spaces at Fig.\ref{lsp}.\footnote{In the diagrams we omit
all $2$-point lines.}
\ecrllr
\def\xx{*=<2mm>[o][F-]{}}
\begin{figure}[h]
$\xymatrix@R=5pt@C=10pt@M=0pt@L=2pt{
\xx &     & \xx & &   \xx\ar@{-}[r] & \xx\ar@{-}[r] & \xx & &    & & \xx\ar@{-}[ddl]\ar@{-}[ddr]     & &   & &   & & \xx \ar@{-}[ddl]\ar@{-}[ddr] & & \\
    &     &     & &                 &               &     & &    & &                                 & &   & &   & &                            & & \\
    &     &     & &                 &               &     & &    & \xx\ar@{-}[ddl] & & \xx\ar@{-}[ddr] &   & &   & \xx\ar@{-}[ddl]\ar@{-}[dr] & & \xx\ar@{-}[dl]\ar@{-}[ddr]& \\
    & \xx &     & &                 & \xx           &     & &    & &                                 & &     & &   & & \xx\ar@{-}[dll]\ar@{-}[drr] & & \\
    &     &     & &                 &               &     & &    \xx & & &                         & \xx   & &   \xx   & & & & \xx \\
}$
\caption{The linear spaces with $\le 7$ points, and $\ge 2$ lines of size $\le 3$, in which
the union of all lines incident to a point coincides with the point set}\label{lsp}
\end{figure}
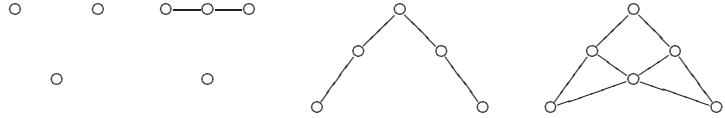
\proof To prove the second statement suppose that $R_{ij}\ne G$ for all $i,j\in I$. Then
(a) two distinct points of~$\G$ are incident to exactly one line, i.e. $\G$ is a linear
space, and (b) the union of all lines incident to a point coincides with~$I$. Without
loss of generality we can assume that $|L|\ge 2$. This implies that each line is
incident at most~$3$ points (for otherwise, any point not in the line is incident to at
least $4$ points in contrast to the first statement). Thus $|I|\le 7$ and the
required statement follows from the list of linear spaces on at most~$9$ points
given in \cite[pp.190-191]{BB93}.\bull

The first linear space at Fig.~\ref{lsp} is known as {\it near-pencil} on $3$ points.

\section{Quasi-thin schemes. Orthogonals}\label{090209g}

A scheme $\X=(\Omega,S)$ is called {\it quasi-thin} if $S=S_1\cup S_2$. In such a scheme
the product of two basic relations is again a basic relation unless both of them are
{\it thick}, i.e. belong to~$S_2$. By \cite[Lemma~4.1]{H01} given a thick relation~$u$
there exists the uniquely determined basic relation $u^\ddp$ such that
$$
uu^*=\{1_\Omega,u^\ddp\}.
$$
This relation is called the {\it orthogonal} of $u$. It is easily seen that any orthogonal
is a (non-reflexive) symmetric relation. 
The following statement was proved in~\cite{HM02a}.

\lmml{260310b}
Given thick relations $u$ and $v$ in the set $S$ we have
\nmrt
\tm{1} $u^\ddp=v^\ddp$ and $u^\ddp\in S_1$ if and only if either $u^*\mult v=2a+2b$ with
$a,b\in S_1$, or $u^*\mult v=2a$ with $a\in S_2$,
\tm{2} $u^\ddp=v^\ddp$ and $u^\ddp\not\in S_1$ if and only if $u^*\mult v=2a+b$ with $a\in S_1$
and $b\in S_2$,
\tm{3} $u^\ddp\ne v^\ddp$ if and only if $u^*\mult v=a+b$ with $a,b\in S_2$.\bull
\enmrt
\elmm

Given $T\subset S_2$ we set $T^\ddp=\{u^\ddp:\ u\in T\}$. Any element from the set $S^\ddp$
is called an orthogonal of the scheme~$\X$.

\thrml{250810a}
Any quasi-thin scheme with at most one orthogonal is schurian and separable.
\ethrm
\proof Let $\X=(\Omega,S)$ be a quasi-thin scheme. If $S^\ddp=\emptyset$, then this
scheme is regular, and hence $1$-regular. Therefore it is schurian and separable by
Theorem~\ref{210109e}. Thus we can assume that $S^\ddp=\{u\}$ for some non-reflexive
basic relation~$u$. Then $u=u^*$ and $v^*v\subset\{1_\Omega,u\}$ for all $v\in S$.
Therefore the set $\{1_\Omega, u\}$ is closed and coincides with the thin
residue~$T$ of the scheme~$\X$.\medskip

Suppose first that $u\in S_2$. Then $u^\ddp=u$ is a thick relation. By statement~(2) of
Lemma~\ref{260310b} this implies that $S_2=S_1u=uS_1$. Therefore $S=S_1T$. Since
$S_1\cap T=\{1_\Omega\}$ and $|T|=2$, it follows that
$$
\X\cong\X_{\alpha S_1}\otimes\X_{\alpha T},\qquad \alpha\in\Omega.
$$
However, the scheme $\X_{\alpha S_1}$ is regular whereas the scheme $\X_{\alpha T}$ is
trivial. Thus both of these scheme are schurian and separable, and we are done
by Theorem~\ref{240810a}.\medskip

Let $u\in S_1$. Then $|\alpha T|=2$ for all $\alpha\in\Omega$. Since every fiber of the
thin residue extension $\X_0=\X_{(T)}$ is of the form $\alpha T$, it follows that
given $\Delta,\Gamma\in\Fib(\X_0)$ the set $\Delta\times\Gamma$ is either a basic relation
of~$\X_0$, or the union of two thin basic relations of $\X_0$. In the latter
case we will write $\Delta\sim\Gamma$. It is easily seen that $\sim$ is an equivalence
relation on the set $\Fib(\X_0)$. Denote by $I$ the set of its classes, and
given $i\in I$ set $\Omega_i$ to be the union of fibers belonging the class~$i$. Then
$$
\X_0=\boxplus_{i\in I}\X_i
$$
where $\X_i=(\X_0)_{\Omega_i}$.
Any summand here is a $1$-regular coherent configuration, and hence is schurian and separable.
By Theorem~\ref{240810a} this implies that so is the coherent configuration~$\X_0$.
Thus the scheme $\X$ is schurian and separable by Theorem~\ref{150610b}.\bull

From \cite[pp.71,72]{W76} it follows that any primitive\footnote{A scheme on $\Omega$ is
called primitive if any equivalence relation on~$\Omega$ that is a union of basic
relations, is $1_\Omega$ or $\Omega\times\Omega$.} quasi-thin scheme is schurian and
separable. Moreover, an inspection of the Hanaki-Miyamoto list~\cite{HM} shows
that any imprimitive quasi-thin scheme of degree at most~$8$ has at most one orthogonal.
Thus by Theorem~\ref{250810a} we have the following statement.

\crllrl{290810a}
Any quasi-thin scheme of degree at most~$8$ is schurian and separable.\bull
\ecrllr

We recall that a scheme is called {\it Kleinian} if its thin residue consists of thin relations
and forms a Klein group with respect to the relation product. In the following statement
these schemes are characterized by means of orthogonals. Below given $u\in S$ we set
$$
S_u=\{v\in S_2:\ v^\ddp=u\}.
$$
Clearly $S_uS_1=S_u$.

\lmml{300810a}
A quasi-thin scheme $\X$ is Kleinian if and only if $S^\ddp\subset S_1$ and either
$|S^\ddp|=2$ or $|S^\ddp|=3$ and the set $\{1_\Omega\}\cup S^\ddp$ is closed. If $\X$ is a
commutative Kleinian scheme, then $|S^\ddp|=3$.
\elmm
\proof The necessity is obvious. To prove the sufficiency without loss of generality
we can assume that $\X$ is a quasi-thin scheme with exactly two thin
orthogonals~$u$ and~$v$. Then given a thick relation~$x\in S_u$, the relation $vx\in S$ is
also thick and
$$
(vx)(vx)^*=v(xx^*)v=v\{1_\Omega,u\}v=\{1_\Omega,vuv\}.
$$
Since $vuv$ is also a basic relation, we conclude that $vuv=(vx)^\ddp\in\{u,v\}$.
Therefore $vuv=u$, i.e.  $u$ and $v$ commute. Thus the
thin residue of $\X$ coincides with the group $\lg u,v\rg=\{1_\Omega,u,v,uv\}$ which
in our case is obviously the Klein group.\medskip

To complete the proof suppose on the contrary that $\X$ is a commutative
Kleinian scheme with exactly two orthogonals~$u$ and~$v$. Then
\qtnl{310810c}
S=S_1\cup S_u\cup S_v\quad\text{and}\quad S_u^*=S_u,\ S^*_v=S_v.
\eqtn
Moreover, given basic relations $x$ and $y$ such that $x^\ddp=y^\ddp$, and any
$z\in xy$ we have $zz^*\subset (xy)(xy)^*=xx^*yy^*\subset\{1_\Omega,x^\ddp\}$.
Since also $S_u=S_uS_1=S_1S_u$ and $S_v=S_vS_1=S_1S_v$, we see that $S_1\cup S_u$
and $S_1\cup S_v$ are closed subsets of~$\X$ the union of which equals~$S$.
Therefore one of them coincides with $S$. A contradiction.~\bull

Let $u\neq v$ be thick basic relations of a quasi-thin scheme~$\X$. We say
that they are {\it adjacent}, $u\approx v$,
if $|u^*v|=2$.  Since $|u^*v|=|v^*u|$, the adjacency relation is symmetric.
Notice that by Lemma~\ref{260310b} the cardinality of $|u^*v|$ is either one or two.
Therefore two relations $u,v\in S_2$ are non-adjacent if and only if $|u^*v|=1$.
The following special statement will be used in the proof of Theorem~\ref{280109e}.

\lmml{290810i}
Let $(\Omega,S)$ be a quasi-thin non-Kleinian scheme of degree $\ge 9$ and with at least
two orthogonals. Suppose that a set $T\subsetneq S_2$ is such that
\qtnl{310810a}
|T^\ddp|\le 2\qquad\text{and}\qquad|T_u|\le 2\quad\text{for all}\quad u\in T^\ddp\cap S_2
\eqtn
where $T_u=T\cap S_u$. Then there exists a relation $t\in S_2\setminus T$ adjacent to each
element of $T$.
\elmm
\proof We observe that if $|T^\ddp|<|S^\ddp|$, then by statement~(3) of Lemma~\ref{260310b}
the required statement holds for any relation $t\in S$ such that
$t^\ddp\in S^\ddp\setminus T^\ddp$. Thus without loss of generality we can assume
that $|T^\ddp|=|S^\ddp|$. Since $|T^\ddp|\le 2$ and $|S^\ddp|\ge 2$, this implies
that $T^\ddp=S^\ddp$ and $|S^\ddp|=2$. If $S^\ddp\subseteq S_2$,
then by statement~(3) of Lemma~\ref{260310b} any two elements of $S_2$
are adjacent, and we are done with arbitrary $t\in S_2\setminus T$. Moreover, taking into
account that the scheme $(\Omega,S)$ is non-Kleinian, we conclude by Lemma~\ref{300810a}
that $S^\ddp\not\subset S_1$. Thus we can assume that
$$
S^\ddp=\{u,v\},\quad u\in S_2,\quad v\in S_1.
$$
If $S_u\setminus T_u\neq\emptyset$, then there exists a relation $t\in S_u\setminus T_u$. By
statement~(3) of Lemma~\ref{260310b} this relation is adjacent to every element of
$S_2\setminus\{t\}$, and we are done. Thus we may assume that $S_u\subseteq T_u$, or
equivalently, $S_u=T_u$. To complete the proof we have to verify that the equality
\qtnl{310810b}
|S_u|= 2
\eqtn
leads to a contradiction. Indeed, let us fix a relation $t\in T_u$. Then $rt=st$ for some
thin relations $r$ and $s$, only if they are equal (here $s^*r\in tt^*=\{1_\Omega,u\}$
and hence $s^*r=1_\Omega$ because $n_u=2$). Moreover, since $u$ is thick, we have $S_1S_u=S_u$.
Therefore from~\eqref{310810b} it follows that $|S_1|=|S_1t|\le|S_1S_u|=|S_u| = 2$.
This implies that
$$
S_1=\{1_\Omega,v\}\quad\text{and}\quad S_u=\{t,vt\}.
$$
Suppose that $u^\ddp=u$. Then $S_u=\{u,vu\}$ and $S_u^*=S_u$. It follows that the set
$Q=S_1\cup S_u$ is closed. However, the set $R=S_1\cup S_v$ is also closed. Thus
the set~$S$ is a union of the closed subsets~$Q$ and~$R$. This implies that one of
them coincides with~$S$ which is impossible. Thus $u^\ddp=v$. Then one can check that
$$
tv=vt\quad\text{and}\quad t^*\in\{t,vt\}.
$$
So $Q=\lg t\rg=S_1\cup S_u\cup\{u\}$ is a closed set and $n_Q=8$. Again the set $S$ is a
union of the closed subsets~$Q$ and~$R$. This implies that $S=Q$ whence it follows that
$|\Omega|=n_S=n_Q=8$. Contradiction.\bull

\section{One-point extension of a quasi-thin scheme}\label{030910b}

In this section we first compute the fibers of a one-point extension of a quasi-thin scheme,
then analyze its basic relations, and finally give a sufficient condition for its
schurity and separability.

\thrml{190410a}
Let $\X=(\Omega,S)$ be a quasi-thin scheme and $\alpha\in\Omega$. Then each fiber of the
coherent configuration $\X_\alpha$ is of the form $\alpha u$, $u\in S$. In
particular,
$$
S_\alpha=\bigcup_{u,v\in S}S_\alpha(u,v)
$$
where $S_\alpha(u,v)=\{a\in S_\alpha:\ a\subset \alpha u\times\alpha v\}$.
\ethrm
\proof The second statement immediately follows from the first one. To prove the
latter let us define an involution $f_\alpha\in\sym(\Omega)$ so that
\qtnl{180410a}
\beta^{f_\alpha}=
\css
\beta, &\text{if $r(\alpha,\beta)\in S_1$,}\\
\beta',&\text{if $r(\alpha,\beta)\in S_2$,}\\
\ecss
\eqtn
where $\beta'$ is defined from the condition $\{\beta,\beta'\}=\alpha r(\alpha,\beta)$.
It was proved in \cite[Lemma~3.5]{HM02b} that $f_\alpha\in\aut(\X)$ for all~$\alpha$.
To complete the proof, let $\Delta$ be a fiber of the coherent configuration $\X_\alpha$.
Then obviously $\Delta\subset\alpha u$ for some $u\in S$. On the other hand, the set
$\alpha u$ is the orbit of the group
$$
\lg f_\alpha\rg\le\aut(\X_\alpha).
$$
Thus $\Delta\supset \alpha u$. Since the converse inclusion is trivial, we conclude
that $\Delta=\alpha u$ and we are done.\bull

The conclusion of Theorem~\ref{190410a} holds for any schurian scheme, and together
with the transitivity of the automorphism group implies the schurity of the scheme
in question. Thus as a consequence of that theorem we obtain the following well-known
statement~\cite{HM02a}.

\crllrl{180410b}
A quasi-thin scheme $\X$ is schurian if and only if the group $\aut(\X)$ is transitive.\bull
\ecrllr

Let $\X=(\Omega,S)$ be a quasi-thin scheme, $\alpha\in\Omega$ and $u,v\in S$. From
Theorem~\ref{190410a} it follows that $1_{\alpha u}\in S_\alpha$. Therefore given
$a\in S_\alpha(u,v)$ the number $|\beta a|=c_{aa^*}^{1_{\alpha u}}$ does not depend
on $\beta\in \alpha u$. Since $|\alpha u|\le 2$ and $|\alpha v|\le 2$, this implies that
\qtnl{200410c}
S_\alpha(u,v)=\{\alpha u\times\alpha v\}\quad\text{or}\quad
S_\alpha(u,v)=\{f_1,f_2\}
\eqtn
where $f_1$ and $f_2$ are the two distinct bijections from $\alpha u$ onto $\alpha v$
(treated as binary relations on $\alpha u\times\alpha v$). Thus the set $S_\alpha(u,v)$
consists of one or two elements, and the latter holds only if $|\alpha u|=|\alpha v|=2$. In this case
the the element of $S_\alpha(u,v)$ other than $a=f_i$ is denoted by $\ov a$.\medskip

Let $u$ and $v$ be basic relations of the scheme~$\X$. Due to~(\ref{150410a})
the intersection numbers $c_{uw}^v$ and $c_{u^*v}^w$ are zero or not simultaneously.
Therefore given $\alpha\in\Omega$ the cardinality of the set
\qtnl{210109f}
S(u,v;\alpha)=\{w_{\alpha u,\alpha v}:\ w\in S,\ c_{uw}^v\ne 0\}
\eqtn
equal the number $|u^*v|$, and hence does not depend on~$\alpha$. The above set
consists of non-empty and pairwise disjoint relations from $S_\alpha(u,v)^\cup$, the
union of which coincides with the set $\alpha u\times\alpha v$. It is easily seen that
$S(u,v;\alpha)^*=S(v,u;\alpha)$.

\lmml{200410e}
Let $(\Omega,S)$ be a quasi-thin scheme with at least two orthogonals,
$\alpha\in\Omega$ and $u,v\in S$. Then $S_\alpha(u,v)\ne S(u,v;\alpha)$ if and only if
$u\not\approx v$. Moreover, in this case $u^\ddp=v^\ddp\in S_1$ and
\qtnl{210410j}
S_\alpha(u,v)=S(u,w;\alpha)\cdot S(w,v;\alpha)
\eqtn
for any $w\in S_2$ with $w^\ddp\ne u^\ddp$.
\elmm
\proof Both $S_\alpha(u,v)$ and $S(u,v;\alpha)$ forms a partition of the set
$\alpha u\times\alpha v$, and the former partition is a refinement of the latter
one. Therefore
\qtnl{290810k}
2\ge |S_\alpha(u,v)|\ge |S(u,v;\alpha)|=|u^*v|.
\eqtn
Thus $S_\alpha(u,v)\ne S(u,v;\alpha)$ if and only if $|S_\alpha(u,v)|=2$ and $|u^*v|=1$.
Due to~(\ref{200410c}) the first equality holds only if $u,v\in S_2$, whereas the second
one means that $u\not\approx v$. This proves the necessity of the first
statement.\medskip

Conversely,
suppose that $u$ and $v$ are non-adjacent elements of~$S_2$. Then $u^\ddp=v^\ddp$ by Lemma~\ref{141208j}.
To complete the proof, let $w\in S_2$ be such that $w^\ddp\ne u^\ddp$. Then
statement~(3) of Lemma~\ref{260310b} implies that $|u^*w|=2$, and hence by~\eqref{290810k}
with $v=w$, we obtain that $|S_\alpha(u,w)|=2$. Similarly, $|S_\alpha(v,w)|=2$. Due
to~(\ref{200410c}) this implies
$$
S_\alpha(u,w)=\{f_1,f_2\}\quad\text{and}\quad S_\alpha(w,v)=\{g_1,g_2\}
$$
where $f_1$ and $f_2$ (resp. $g_1$ and $g_2$) are the two bijections from $\alpha u$
to $\alpha w$ (resp. from $\alpha w$ to $\alpha v$). Since obviously
$S(u,w;\alpha)\cdot S(w,v;\alpha)\subset S_\alpha(u,v)^\cup$ and
$$
f_i\cdot g_j=f_ig_j,\qquad i,j=1,2,
$$
this implies that $f_ig_j\in S_\alpha(u,v)$. This proves equality~(\ref{210410j}),
and also the sufficiency of the first statement.\bull

\crllrl{240210a}
Let $\X$ be a quasi-thin scheme on $\Omega$ with at least two orthogonals. Then
given a point $\alpha\in\Omega$ the coherent configuration $\X_\alpha$ is $1$-regular.
\ecrllr
\proof Denote by $S$ the set of basic relations of the scheme~$\X$. Let us verify that any point
$\beta\in\alpha S_2$ is regular (see Subsection~\ref{250210c}). To do this let $a\in S_\alpha$ be such that
$\beta a\ne\emptyset$. Then by Theorem~\ref{190410a} there exist relations $u\in S_2$
and $v\in S$ such that
$$
\beta\in \alpha u\quad\text{and}\quad a\in S_\alpha(u,v).
$$
Without loss of generality we can assume that $v\in S_2$ (otherwise $|\beta a|=|\alpha v|=1$
and we are done). Then due to~\eqref{200410c} it suffices to verify that $|S_\alpha(u,v)|=2$.
However, this is true by Lemma~\ref{200410e}:
if $S_\alpha(u,v)=S(u,v;\alpha)$, then $|S_\alpha(u,v)|=|u^*v|=2$, whereas
if $S_\alpha(u,v)\ne S(u,v;\alpha)$, then $2\ge |S_\alpha(u,v)|>|u^*v|=1$.\bull

The conclusion of Corollary~\ref{240210a} is not true when $|S^\ddp|=1$. Indeed,
denote by $\X$ the scheme of the wreath product of two regular schemes of degrees~$2$
and~$n\ge 3$. Then $\X$ is a quasi-thin scheme of degree~$2n$ with exactly one orthogonal.
On the other hand, any
point extension of $\X$ is the coherent configuration of the elementary abelian
group of order~$2^{n-1}$ with two fixed points and $n-1\ge 2$ orbits of cardinality~$2$.
It follows that the point extension of~$\X$ has no regular points, and hence is not
$1$-regular.

\thrml{250210r}
Let $\X$ be a quasi-thin scheme with at least two orthogonals. Suppose
that any algebraic isomorphism $\varphi$  from $\X$ to another scheme $\X'$
has one point extension $\varphi_{\alpha,\alpha'}:\X_{\alpha}\to\X'_{\alpha'}$ for any
pair of points $\alpha\in\Omega$ and $\alpha'\in\Omega'$. Then the scheme $\X$ is schurian
and separable.
\ethrm
\proof
By Corollary~\ref{240210a} the coherent configuration $\X_\alpha$ is $1$-regular.
Together with Theorem~\ref{210109e} this implies that the set
$\iso(\X^{}_{\alpha^{}},\X'_{\alpha'},\varphi_{\alpha^{},\alpha'})$ is not empty. Since
$$
\iso(\X^{}_{\alpha^{}},\X'_{\alpha'},\varphi_{\alpha^{},\alpha'})\subset\iso(\X,\X',\varphi),
$$
the set $\iso(\X,\X',\varphi)$ is also not empty. Thus the scheme $\X$ is separable.\medskip

To prove schurity of $\X$ take $\alpha,\alpha'\in\Omega$. Then by the theorem hypothesis
the trivial algebraic isomorphism $\id_S:\X\to\X$ has the $(\alpha,\alpha')$-extension,
say $\varphi_{\alpha,\alpha'}$. Since the coherent configuration
$\X_\alpha$ is $1$-regular (Corollary~\ref{240210a}), from Theorem~\ref{210109e} it
follows that there exists an isomorphism
$$
f_{\alpha,\alpha'}\in\iso(\X_{\alpha^{}},\X_{\alpha'},\varphi_{\alpha,\alpha'}).
$$
By the definition of $\varphi_{\alpha,\alpha'}$ this isomorphism takes $\alpha$ to $\alpha'$
and preserves every basic relation of~$\X$. Therefore $f_{\alpha,\alpha'}\in\aut(\X)$. Since $\alpha$ and $\alpha'$ are arbitrary
points of $\Omega$, this means that the group $\aut(\X)$ is transitive. Thus schurity
of~$\X$ follows from statement~(2) of Theorem~\ref{180410b}.~\bull

\section{Triangles in a quasi-thin scheme}\label{180610a}
Let $\X=(\Omega,S)$ be a quasi-thin scheme. A $3$-subset $T$ of $S_2$ is called a
{\it triangle} (in~$\X$) if any two distinct elements of $T$ are adjacent.  From Lemma~\ref{200410e} it
follows that any $3$-set $T\subset S_2$ with $|T^\ddp|=3$ is a triangle. The following
statement which is an immediate consequence of~(\ref{200410c}), shows that any triangle induces a regular
coherent configuration with three fibers of size~$2$ on the neighborhood of each point.

\lmml{190510a}
Let $\{u,v,w\}$ be a triangle in the scheme~$\X$ and $\alpha\in\Omega$. Then
$x\cdot y=\ov x\cdot\ov y$ and $\ov x\cdot y=x\cdot\ov y=\ov{x\cdot y}$
for all $x\in S_\alpha(u,w)$ and $y\in S_\alpha(w,v)$.\bull
\elmm

We say that a triangle $\{u,v,w\}$ is {\it exceptional} if $u^\ddp\cdot v^\ddp\cdot w^\ddp=1_\Omega$.
Notice that in this case $u^\perp,v^\perp,w^\perp$ are pairwise distinct thin elements which
together with the identity form a thin closed subset of $S$ isomorphic to a Klein group. Conversely,
if $T$ is a triangle for which the set $1_\Omega\cup T^\ddp$ is a Klein subgroup of~$S_1$,
then obviously the triangle $T$ is exceptional. Thus  we come to the following statement.

\thrml{190510k}
A triangle $T$ is exceptional if and only if the set $1_\Omega\cup T^\ddp$ is a
Klein subgroup of the group~$S_1$.\bull
\ethrm

 The following theorem
provides the key property of non-exceptional triangles that can be used in Section~\ref{310309b}.

\thrml{regular}
Let $T = \{u,v,w\}$ be a non-exceptional triangle. Then there exist
relations $a\in u^*w$ and $b\in w^*v$ for which $|ab\cap u^*v|=1$.
\ethrm
\proof Pick arbitrary $a\in u^*w$ and $b\in w^*v$. Then $w\in ua\cap vb^*$. Therefore
$|ua\cap vb^*|\geq 1$, and hence $|ab\cap u^*v|\geq 1$. If one of the sets $u^*w$ or $w^*v$
contains a thin element, then we can choose $a$ or $b$ to be thin. But then
$|ab|=1$ and we are done in this case. Thus for the rest of the proof we can
assume that $u^*w, w^*v\subseteq S_2$. Notice that by Lemma~\ref{260310b} this implies that
$w^\ddp\neq u^\ddp$ and $w^\ddp\neq v^\ddp$.\medskip

Assume now, towards a contradiction, that $|ab\cap u^*v|\geq 2$ for all
$a\in u^* w$ and $b\in w^*v$. Then $2\le |ab\cap u^*v|\le |u^*v|=2$
for all $a$ and $b$. So $|ab\cap u^*v|=2$. Therefore
\qtnl{040610a}
ab=u^*v,\qquad a\in u^* w,\ b\in w^*v.
\eqtn
If now $u^\perp=v^\perp$, then by Lemma~\ref{260310b} the set $u^*v$ contains a thin
element, say~$t$. This implies that $t\in ab$ for every $a\in u^*v$ and a fixed $b\in w^*v$.
But then $u^*v$ consists of the unique element $a=tb^*$. Contradiction. Thus
$u^\perp\ne v^\perp$, and hence $u^*v\subset S_2$. Due to (\ref{040610a})
this shows that $ab$ and $u^*v$ are equal as multisets. So
\qtnl{190510f}
2u^*\mult v+\frac{2}{n_{w^\ddp}}\mult u^*\mult w^\ddp v=u^*(ww^*)v=(u^*w)(w^*v)=\sum_{a\in u^*w,b\in w^*v} ab=4u^*v,
\eqtn
implying $u^* w^\perp v=n_{w^\perp} u^*v$. Suppose first that $n_{w^\ddp}=2$. Then
$u^*w^\ddp v=2u^*v$. Therefore one can find
points $\alpha_i,\beta_1,\beta_2\in\Omega$ where $i=0,\ldots,3$ and $\beta_1\ne\beta_2$,
to have the configuration at Fig.\ref{fg}.
\begin{figure}[h]
$\xymatrix@R=30pt@C=15pt@M=0pt@L=2pt{
& & & \VRT{\beta_1} \ar[drrr]^*{v} & & & \\
\VRT{\alpha_0} \ar[urrr]^*{u^*} \ar[drrr]_*{u^*} \ar[rr]^*{u^*} & & \VRT{\alpha_1} \ar[rr]_*{w^\ddp} & &
\VRT{\alpha_2} \ar[rr]^*{v} & & \VRT{\alpha_3}\\
& & & \VRT{\beta_2} \ar[urrr]_*{v} & & & \\
}$
\caption{}\label{fg}
\end{figure}
However, $n_{u^*}=2$. So either $\alpha_1=\beta_1$
or $\alpha_1=\beta_2$. But then in any case $r(\alpha_1,\alpha_3)$ is contained
in the set $vv^*\cap w^\ddp$ which is empty because $v^\ddp\ne w^\ddp$. Contradiction.
Thus $n_{w^\ddp}=1, w^\ddp\in S_1$ and hence by~(\ref{190510f}) we have $u^*w^\ddp v=u^* v$.
Therefore  $\lg u^*w^\ddp v,u^* v\rg=\lg u^* v,u^* v\rg=4$ which together with~\eqref{associative} implies
that
$$
4=\lg uu^* w^\ddp,vv^*\rg=
\lg 2w^\ddp+k\,u^\ddp w^\ddp,2\,1_\Omega+m\, v^\ddp\rg=
km\lg u^\ddp w^\ddp, v^\ddp\rg
$$
where $k,m\in\{1,2\}$. In particular, the right-hand side of the above equality
is non-zero. However, since $u^\ddp\ne w^\ddp$ and $w^\ddp\in S_1$, the latter is possible
only if $u^\ddp w^\ddp=v^\ddp$ and $k=m=2$ meaning that $T$ is exceptional. Contradiction.\bull

\section{One point extension of an algebraic isomorphism}\label{310309b}

\sbsn\label{010510b}
In this section we prove the following theorem which is the key ingredient in
the proof of our main result.

\thrml{280109e}
Let $\X$ be a non-Kleinian quasi-thin scheme of degree $\ge 9$. Suppose that it
has at least two orthogonals. Then
any algebraic isomorphism $\varphi$  from $\X$ to another scheme $\X'$
has a one point extension $\varphi_{\alpha,\alpha'}:\X_{\alpha}\to\X'_{\alpha'}$
for any pair of points $\alpha\in\Omega$ and $\alpha'\in\Omega'$.
\ethrm
\proof Let $\X=(\Omega,S)$ and let $\varphi:u\mapsto u'$ be an algebraic isomorphism from
$\X$ to a scheme $\X'=(\Omega',S')$. Then it is easily seen that $\X'$ is quasi-thin
and
\qtnl{200410a}
(u^\ddp)'=(u')^\ddp,\quad |u^*v|=|(u')^*v'|,\qquad u,v\in S.
\eqtn
Let us fix points $\alpha\in\Omega$ and $\alpha'\in\Omega'$. In the following two subsections
we will construct a bijection
\qtnl{210410i}
\varphi':S^{}_{\alpha^{}}\to S'_{\alpha'},\ a\mapsto a',
\eqtn
such that $(S^{}_{\alpha^{}}(u,v))'=S'_{\alpha'}(u',v')$ for all $u,v\in S$. In
Subsection~\ref{020510a} we define $\varphi'$ on the union of all sets $S_\alpha(u,v)$
with $u\approx v$; in Subsection~\ref{020510b} we extend the obtained
mapping on the set $S_\alpha$. In Subsection~\ref{020510c} it will be proven that
$\varphi'$ is the $(\alpha,\alpha')$-extension of~$\varphi$.

\sbsn\label{020510a}
Let $u,v\in S$ be such that $S_\alpha(u,v)=S(u,v;\alpha)$. Then any
$a\in S_\alpha(u,v)$ is of the form $a=w_{\alpha u,\alpha v}$ for some $w\in S$.
Set
\qtnl{210410z}
a'=w'_{\alpha' u',\alpha' v'}.
\eqtn
Then $a'\in S'(u',v';\alpha')$. Due to~(\ref{200410a}) we have $|u^*v|=|(u')^*v'|$.
Therefore the mapping $a\mapsto a'$ is a bijection. It is easily seen that if
$S_\alpha(u,v)=\{a,\ov a\}$, then $(\ov a)'=\ov{a'}$.

\lmml{210410g}
Given a triangle $\{u,v,w\}\subset S_2$ and $a\in S_\alpha(u,v)$,
$b\in S_\alpha(v,w)$, $c\in S_\alpha(u,w)$ we have
$$
a\cdot b=c\ \Rightarrow\ a'\cdot b'=c'.
$$
\elmm
\proof Assume first that the triangle $\{u,v,w\}$ is non-exceptional. Then
by Theorem~\ref{regular} there exist
relations $x\in u^*w$, $y\in w^*v$ and $z\in u^*v$ such that
$xy\cap u^*v =\{z\}$. This implies that $x'y'\cap {u'}^*v'=\{z'\}$, and
$a_1\cdot b_1=c_1$ and $a'_1\cdot b'_1=c'_1$ where
$$
a_1=x_{\alpha u,\alpha w},\quad b_1=y_{\alpha w,\alpha v},\quad
c_1=z_{\alpha u,\alpha v}
$$
Now let $a\in S_\alpha(u,v)$, $b\in S_\alpha(v,w)$, $c\in S_\alpha(u,w)$ be such that
$a\cdot b=c$. There the pair $(a,b)$ is one of the following: $(a_1,b_1)$,
$(\ov a_1,b_1)$, $(a_1,\ov b_1)$ or $(\ov a_1,\ov b_1)$. In the first case the statement is
clear. In the second one we have $a'=(\ov a_1)'=\ov{a'_1}$. By Lemma~\ref{190510a} this
implies that $c=a\cdot b=\ov c_1$. Thus
and
$$
a'\cdot b'=\ov{a'_1}\cdot b'_1=\ov{a'_1\cdot b'_1}=\ov{c'_1}=\ov{c_1}'=c'.
$$
The remaining two cases are considered in a similar manner.\medskip

Let now $T=\{u,v,w\}$ be an exceptional triangle. Since $\X$ is a non-Kleinian scheme,
this implies that $|S^\ddp|\geq 4$. Therefore there exists a thick basic relation~$t$
such that $t^\ddp\not\in T^\ddp$. Then each set $\{t,x,y\}$ where $x$ and $y$ are distinct
elements of $T$, is a non-exceptional triangle (otherwise $t^\ddp=x^\ddp\cdot y^\ddp\in T^\ddp$).
Therefore there exist $x\in S_\alpha(t,u)$, $y\in S_\alpha(t,v)$ and $z\in S_\alpha(t,w)$
for which $x^*\cdot y=a$, $y^*\cdot z=b$ and $x^*\cdot z=c$. Since the corresponding
triangles are non-exceptional, from the first part of the proof it follows that
$$
x'^*\cdot y'=a',\quad y'^*\cdot z'=b',\quad x'^*\cdot z'=c'.
$$
Therefore $a'b'=(x'^*\cdot y')\cdot(y'^*\cdot z')=x'^*\cdot z'=c'$. The lemma is proved.\bull

\sbsn\label{020510b}
Let $u,v\in S$ be such that $S_\alpha(u,v)\ne S(u,v;\alpha)$. Then by Lemma~\ref{200410e}
the relations $u$ and $v$ are non-adjacent and $u^\ddp=v^\ddp$. Since the scheme~$\X$ has at least two
orthogonals, one can find a relation $w=w_{u,v}$ for which $w^\ddp\ne u^\ddp$. Then by
Lemma~\ref{200410e} we have
\qtnl{210410x}
S_\alpha(u,v)=S(u,w;\alpha)\cdot S(w,v;\alpha).
\eqtn
Here $u\approx w$ and $w\approx v$, and $S_\alpha(u,w)=S(u,w;\alpha)$ and
$S_\alpha(w,v)=S(w,v;\alpha)$. Therefore one can consider two bijections
$$
S_\alpha(u,w)\to S'_{\alpha'}(u',w'),\ b\mapsto b',\qquad
S_\alpha(w,v)\to S'_{\alpha'}(w',v'),\ c\mapsto c'
$$
that were defined in Subsection~\ref{020510a}. Now for a fixed $b\in S_\alpha(u,w)$
and for any $a\in S_\alpha(u,v)$ there exists a uniquely determined $c\in S_\alpha(w,v)$
such that $a=b\cdot c$. Set
\qtnl{210410w}
a'=b'\cdot c'.
\eqtn
Then $a'\in S'(u',w';\alpha')\cdot S'(w',v';\alpha')=S'_{\alpha'}(u',v')$ and the
mapping $a\mapsto a'$ is a required bijection.

\lmml{260109o}
In the above notation set $w_1=w$, $b_1=b$ and $c_1=c$.  Then given $w_2\in S_2$
with $w_2^\ddp\ne u^\ddp$, $b_2\in S_\alpha(u,w_2)$ and $c_2\in S_\alpha(w_2,u)$ we have
$$
b_1\cdot c_1=b_2\cdot c_2\ \Rightarrow\ b'_1\cdot c'_1=b'_2\cdot c'_2.
$$
\elmm
\proof Without loss of generality we can assume that $w_1\ne w_2$. In the above
assumptions the element $u^\ddp=v^\ddp$ is not equal neither to $w_1^\ddp$ nor to $w_2^\ddp$.
Therefore $u\approx w_1$, $w_1\approx v$, $v\approx w_2$ and $w_2\approx u$.\medskip

Suppose that $w_1\approx w_2$. It follows from the equality $b_1\cdot c_1=b_2\cdot c_2$
that the relation $d:=b_2^*\cdot b_1=c_2\cdot c_1^*$ belongs to the set $S_\alpha(w_2,w_1)$.
 Then $b_1=b_2\cdot d$ and $c_1=d^*\cdot c_2$. Since both $\{w_1,w_2,u\}$ and $\{w_1,w_2,v\}$
are triangles, Lemma~\ref{210410g} implies
$$
b'_1\cdot c'_1=(b_2\cdot d)'\cdot (d^*\cdot c_2)'=b_2'\cdot d'\cdot (d')^*\cdot c_2'=
b'_2\cdot c'_2
$$
and we are done.\medskip

Let now $w_1\not\approx w_2$. Then $w_1^\ddp=w_2^\ddp$ and this element is thin
(Lemma~\ref{260310b}). By Lemma~\ref{290810i} applied to $T=\{u,v,w_1,w_2\}$ there exists a
relation $t\in S_2$ such that any set $\{t,x,y\}$ with $x\in\{u,v\}$ and $y\in\{w_1,w_2\}$
is a triangle. Pick an arbitrary $a_1\in S_\alpha(u,t)$. Then we have
$$
d_1:=a_1^*\cdot b_1\in S_\alpha(t,w_1),\qquad
d_2:=a_1^*\cdot b_2\in S_\alpha(t,w_2),\qquad
a_2:=d_1\cdot c_1\in S_\alpha(t,v).
$$
Therefore $a_1\cdot a_2=(a_1\cdot d_1)\cdot (d^*_1\cdot a_2)=b_1\cdot c_1=b_2\cdot c_2=a_1\cdot d_2\cdot c_2$.
This implies that $c_2=d_2^*\cdot a_2$ (Figure~\ref{f5}).
\begin{figure}[h]
\grphp{
 & & & \VRT{w_1} \ar[dddrrr]^*{c_1} & & & \\
 & & & & & & \\
 & & & & & & \\
 \VRT{u}\ar[rrruuu]^*{b_1}\ar[rrr]_*{a_1}\ar[rrrddd]_*{b_2} & & &
 \VRT{t}\ar[uuu]_*{d_1}\ar[rrr]_*{a_2}\ar[ddd]^*{d_2} & & &
 \VRT{v} \\
 & & & & & & \\
 & & & & & &\\
 & & & \VRT{w_2} \ar[uuurrr]_*{c_2} & & & \\
}
\caption{}\label{f5}
\end{figure}
Thus by Lemma~\ref{210410g} we conclude that
$$
b'_1\cdot c'_1=(a_1\cdot d_1)'\cdot (d^*_1\cdot a_2)'=
a'_1\cdot (d'_1\cdot (d'_1)^*)\cdot a'_2=a'_1\cdot a'_2=
$$
$$
(b_2\cdot d_2^*)'\cdot (d_2\cdot c_2)'=
b'_2\cdot (d'_2\cdot (d^*_2)')\cdot c'_2=b'_2\cdot c'_2
$$
which completes the proof.\bull\medskip

\sbsn\label{020510c}
For the bijection $\varphi'$ defined in Subsections~\ref{020510a} and~\ref{020510b}
the conditions~(\ref{290309z}) are obviously satisfied. Thus to check that $\varphi'$
is the $(\alpha,\alpha')$-extension of $\varphi$ it suffices to verify only that
\qtnl{300309a}
(b\cdot c)'=b'\cdot c',\qquad b,c\in S_\alpha,\ b\cdot c\ne\emptyset.
\eqtn
Let $b,c\in S_\alpha$ be such that $b\cdot c\ne\emptyset$. Then there exist $u,v,w\in S$ such that
$b\subset \alpha u\times\alpha v$ and $c\subset \alpha v\times\alpha w$. Without loss of
generality we can assume that $u,v,w\in S_2$. If the cardinality of the set
$T=\{u,v,w\}$ is less than three, then at least one of the relations $b, c, b\cdot c$
is an in-fiber relation and we are done by Lemma~\ref{190510a}. So we may assume that
$|T|=3$.\medskip

First we notice that~\eqref{300309a} is correct when $u\approx v$ and $v\approx w$. Indeed,
if in addition $u$ and $w$ are adjacent, then $T$ is a triangle and we are done by
Lemma~\ref{210410g}; otherwise $u$ and $w$ are non-adjacent and we are done
by Lemma~\ref{260109o}. Thus we can assume that at least one of the pairs $\{u,v\}$,
$\{v,w\}$ is non-adjacent. By Lemma~\ref{260310b} this implies that
$$
|T^\ddp|\leq 2.
$$
Now Lemma~\ref{290810i} applied to $T=\{u,v,w\}$ implies that there exists a relation $t\in S_2\setminus T$ such
that any set $\{t,x,y\}$ with $x,y\in T$, is a triangle. As we have shown before, the condition~\eqref{300309a} holds for any
$b\in S_\alpha(x,t)$ and $c\in S_\alpha(t,y)$. On the other hand, by Corollary~\ref{240210a}
the coherent configuration $S_\alpha$ is $1$-regular. Thus there exist relations
$a_1\in S_\alpha(u,t)$, $a_2\in S_\alpha(t,v)$ and $a_3\in S_\alpha(t,w)$
such that $b=a_1\cdot a_2$ and $c=a_2^*\cdot a_3$ (see Fig.~\ref{f6}).
\begin{figure}[h]
\grphp{
 & & & \VRT{t} \ar[ddddrrr]^*{a_3}\ar[dddd]^*{a_2} & & & \\
 & & & & & & \\
 & & & & & & \\
 & & & & & & \\
\VRT{u}\ar[rrruuuu]^*{a_1}\ar[rrr]_*{b} & & &
\VRT{v}\ar[rrr]_*{c} & & &
\VRT{w} \\
}
\caption{}\label{f6}
\end{figure}
Since
$|u^*t| = |t^*v|=2$, the condition~\eqref{300309a} holds for $a_1$ and $a_2$
implying $b'=a'_1\cdot a'_2$. Analogously,  $c'=(a_2^*)'\cdot a'_3$
and $(a_1\cdot a_3)'=a'_1\cdot a'_3$. Therefore
$$
(b\cdot c)'=(a_1\cdot a_2\cdot a_2^*\cdot a_3)'=(a_1\cdot a_3)'=
a'_1\cdot a'_3=
$$
$$
a'_1\cdot a'_2\cdot (a_2^*)'\cdot a'_3=
(a_1\cdot a_2)'\cdot (a_2^*\cdot a_3)'=b'\cdot c'
$$
which completes the proof of~(\ref{300309a}). Theorem~\ref{280109e} is proven.\bull

\section{Proofs of the main results}\label{170610h}

In this section $\X=(\Omega,S)$ denotes a quasi-thin scheme and we write $1$ instead of $1_\Omega$.

\sbsnt{Proof of Theorem~\ref{210210a}}\label{030910a}
Suppose first that $\X$ is non-Kleinian. We have to prove that it
is schurian and separable. By Corollary~\ref{290810a} and Theorem~\ref{250810a}
we can assume that $\X$ is of degree $\ge 9$ and has at least two orthogonals.
Then by Theorem~\ref{280109e} any algebraic isomorphism from $\X$ to another scheme has
one point extension at every pair of points. Thus the scheme $\X$ is schurian and separable
by Theorem~\ref{250210r}.\medskip

Suppose that the scheme $\X$ is Kleinian. Denote by $T$
its thin residue. Then the thin residue extension $\X_0=\X_{(T)}$ is a Klein configuration.
For this configuration we keep the notation of Section~\ref{010710s} with $G=T$ and the
group isomorphisms~\eqref{240610b} taking $g\in T$ to $g^{}_i=g^{}_{\Omega_i,\Omega_i}$,
$i\in I$; in particular, $|I|=|\Omega/T|$. Let
$$
\Psi=\{\psi_f:\ f\in\aut(\X_{\Omega/T})\}
$$
be the group of algebraic automorphisms $\psi_f$ of the scheme~$\X$ defined in statement~(i)
in the proof of Theorem~\ref{150610b}. Then $\Psi$ acts regularly on the set $\Fib(\X_0)$,
and hence on the set $I$ so that $i^\psi=j$ if and only if $(\Omega_i)^\psi=\Omega_j$. By the choice of
isomorphisms~\eqref{240610b} for any $i,j\in I$ and $s\in S$ we have
$$
(s_{ij}\cdot g_j)^\psi=s_{ij}^\psi\cdot(g_j)^\psi=
s_{i^\psi j^\psi}\cdot g_{j^\psi},\qquad\ g\in G,\ \psi\in\Psi,
$$
where $s_{ij}=s_{\Omega_i,\Omega_j}$. This implies that $R_{i^{}j^{}}=R_{i^\psi,j^\psi}$
for all $i,j$ where $R_{ij}$ is the group defined in \eqref{020710a} for $\X=\X_0$. Thus
\qtnl{020710b}
R(\X_0,G)^\Psi=R(\X_0,G).
\eqtn
We note that no entry of this matrix equal to $G$. Indeed, otherwise from the condition~(K1)
it follows that $\X_0$ contains a basic relation $s=\Omega_i\times\Omega_j$
for some $i,j\in I$. However, then the basic relation of the scheme $\X$ that contains~$s$
has valency~$\ge 4$ which is impossible because the scheme $\X$ is quasi-thin. Thus
the hypothesis of Corollary~\ref{010710a} is satisfied.

\lmml{010710u}
The isomorphism type of the linear space $\G_J=\G((\X_0)_J)$ does not depend on the transversal
$J\in\J(\X_0,G)$. Moreover, $\G_J$ is isomorphic to either near-pencil on~$3$ points
or a projective or affine plane of order~$2$.
\elmm
\proof Let $J$ be a transversal of the partition of the set~$I$ in the classes of the
equivalence relation~$\sim$. Then the linear space $\G_J$ has
at least two lines. Indeed, suppose on the contrary that $L_i(H)=I$ for some element $i\in I$
and a group $H\le G$ of order~$2$. Then
$$
R_{ji}=H,\qquad j\in I\setminus\{i\},
$$
where $R=R(\X_0,G)$. Due to~\eqref{020710b} this implies that any non-diagonal entry of the
matrix~$R$ equals to~$H$. Therefore the scheme $\X$ has a unique orthogonal $g_i^\Psi$
where $g$ is the element of~$H$ of order~$2$. However, this is impossible because
$|S^\ddp|\ge 2$. Thus by Corollary~\ref{010710a} the linear space $\G_J$ is a projective or
affine plane of order~$2$, or $\G_J$ is one of the four linear spaces at Fig.\ref{lsp}.\medskip

Given nonnegative integers $d,e$ denote by $M_{d,e}$ the set of all pairs
$(i,J)\in I\times\J(\X_0,G)$ such that $i\in J$ and the linear space $\G_J$ contains
exactly $d$ (resp. $e$) lines of size~$2$ (resp. of size~$3$) that are incident to~$i$.
Suppose that $(i,J)\in M_{d,e}$. Then from condition \eqref{300610z} it follows that
$\{i\}\times\J_i(\X_0,G)\subset M_{d,e}$ where $\J_i(\X_0,G)$ is the set of all
transversals $J\in \J(\X_0,G)$ containing~$i$. Due to~\eqref{020710b} this implies that
$$
\bigcup_{i\in I}\{i\}\times\J_i(\X_0,G)\subset M_{d,e}.
$$
Thus for any $J\in\J(\X_0,G)$ the linear space $\G_J$ contains exactly $d$ (resp. $e$)
lines of size~$2$ (resp. of size~$3$) through any point. However, this is possible
only if $d=e$ and this number is~$2$ or~$3$. In the former case $\G_J$ is the first linear
space at Fig.\ref{lsp} or an affine plane of order~$2$, whereas in the latter case
$\G_J$ is a projective plane of order~$2$. Since all these geometries have distinct
number of points, we are done.\bull

Depending on the isomorphism type of linear spaces $\G_J$ we will say that the scheme~$\X$
is a scheme over near-pencil, affine plane or projective plane. It should be noted that
the number of points in~$\G_J$ coincides with the number~$|J|=|\Omega|/|S_1|$ which
was called the index of~$\X$ in the introduction.

\thrml{020710u}
Any quasi-thin Klein scheme $\X$ is of index $3,4$ or $7$; in these cases $\X$ is
a scheme over near-pencil, affine plane or projective plane. Moreover, in the former
case~$\X$ is schurian and separable, whereas in the latter case~$\X$ is not
commutative.
\ethrm
\proof The first statement immediately follows from Lemma~\ref{010710u}. A straightforward
computation shows that any Klein configuration on~$12$ points is schurian and separable.
Therefore the schurity and separability of a Kleinian quasi-thin scheme over near-pencil
follows from Theorem~\ref{150610b} and Lemma~\ref{290610k}. To prove the last
statement we observe that the commutativity of the scheme~$\X$ implies that the matrix
$R(\X_0,G)$ is symmetric. Therefore either a linear space $\G_J$ contains exactly one line
or it has two disjoint lines. Since both of these possibilities are impossible when
$\G_J$ is a projective plane, we are done.\bull

By Theorem~\ref{020710u} to complete the proof we have to verify that given $i\in\{4,7\}$
there exist infinitely many both non-schurian and non-separable Kleinian schemes of
index~$i$. To do this we denote by $\X_{16}$ and $\X'_{16}$ (resp. $\X_{28}$ and $\X'_{28}$)
the schemes \#173 and \#172 (resp. \#176 and \#175) from the
Hanaki-Miyamoto list~\cite{HM} of association schemes of degree~$16$ (resp. of
degree~$28$). A straightforward computation shows that:
\nmrt
\tm{1} all the schemes $\X^{}_{16}$, $\X'_{16}$, $\X^{}_{28}$ and $\X'_{28}$ are quasi-thin
and Kleinian; the former two are of index~$4$ whereas the latter two are of index~$7$,
\tm{2} the schemes $\X'_{16}$ and $\X'_{28}$ are schurian whereas the schemes $\X^{}_{16}$
and $\X^{}_{28}$ are non-schurian,
\tm{3} the scheme $\X^{}_{16}$ is algebraically isomorphic to the scheme $\X'_{16}$,
the scheme $\X^{}_{28}$ is algebraically isomorphic to the schemed $\X'_{28}$.
\enmrt
Thus $\X_{16}$ and $\X_{32}$ are non-schurian and non-separable quasi-thin Klein schemes
of indices~$4$ and~$7$. Let $\X$ be one of this scheme and let $\Y$ be an arbitrary regular
scheme. Then obviously $\X\otimes\Y$ is a quasi-thin Klein scheme of the same index as~$\X$.
Since by Theorem~\ref{240810a} the scheme $\X\otimes\Y$ is also non-schurian and non-separable,
we are done.\bull

\sbsnt{Proof of Theorem~\ref{180510a}.}\label{030910c}
Suppose that the scheme $\X$ is commutative. To prove that it is schurian by
Theorem~\ref{210210a} we can assume that this scheme is Kleinian. Then from
Lemma~\ref{300810a} it follows that
$$
S^\ddp=\{a,b,c\}
$$
where $a,b,c\in S_1$ with $a^2=b^2=c^2=1$ and $ab=c$, and hence
\qtnl{cond}
S=S_1\cup S_a \cup S_b \cup S_c.
\eqtn
Moreover, given $e\in\{a,b,c\}$ we have $S_e=S_1S_e=S_eS_1$. Choose elements
$x\in S_a$, $y\in S_b$, $z\in S_c$ so that $1\in xyz$. Then
\qtnl{xyz}
\begin{array}{rcl}
xy & = & z^* + z^* a  = z^* + z^* b,\\
yz & = & x^* + x^* b  = x^* + x^* c,\\
zx & = & y^* + y^* c  = y^* + y^* a.
\end{array}
\eqtn
It follows from $x^*\in S_a,y^*\in S_b,z^*\in S_c$ that there exist $u,v,w\in S_1$ such that
\qtnl{star}
x^*=xu^*,\qquad y^*=yv^*,\qquad z^*=zw^*.
\eqtn
Since $x = xa, y=yb, z=zc$, there is a certain freedom in a choice of $u,v,w$.
More precisely, we can always replace (if necessary) $u$ by $ua$, $v$ by $vb$ and $w$ by $wc$.
All these replacements could be done independently.\medskip

Applying ${}*$ to the first row of \eqref{xyz} we obtain that $x^*y^*=z+za$. By \eqref{star}
this implies that
$$
z^* + z^* a=xy=z^*(uvw)+z^*a(uvw).
$$
Therefore $uvw\in\{1,a,b,c\}$.
If $uvw = a$, then by replacing $u$ by $ua$ we obtain that $uvw=1$ (the same could be done
in the cases $uvw=b$ or $uvw = c$). Thus in what follows we can assume that $uvw=1$.\medskip

Let $\alpha,\beta,\gamma\in\Omega$ be such that $(\alpha,\beta)\in x$, $(\beta,\gamma)\in y$ and
$(\gamma,\alpha)\in z$. Since $c_{xy}^{z^*}=c_{yx}^{z^*}=1$, there exists a unique point
$\delta\in\Omega$ such that $(\alpha,\delta)\in x$ and $(\delta,\gamma)\in y$.
The pair $(\delta,\beta)$ belongs to the relation
$$
yx^*=yxu^*=z^*u^*+z^*u^*b=zw^*u^*+zw^*u^*b=zv+zvb,
$$
and hence belongs to either $zv$ or $zvb$. The latter case may be reduced
to the first one by the replacement $u\leftrightarrow ua$, $v\leftrightarrow vb$,
$w\leftrightarrow wc$. Notice that this replacement keeps invariant the relations $1\in xyz$
and $uvw=1$. Thus
we may assume that $(\delta,\beta)\in zv$. This yields us the picture on Fig.~\ref{f2}.
\begin{figure}[h]
$\xymatrix
{
\VRT{\beta}\ar[rr]^x &  & \VRT{\gamma}\ar[lld]^z\\
\VRT{\alpha}\ar[u]^y\ar[rr]_x & & \VRT{\delta}\ar[u]_y\ar[ull]_{zv}\\
}$
\caption{}\label{f2}
\end{figure}
Thus $r(\lambda,\mu)\in S_2$ for all distinct elements $\lambda$, $\mu$ in the set
$\Lambda=\{\alpha,\beta,\gamma,\delta\}$. Therefore due to (\ref{cond}) the set
$\Omega$ is a disjoint union of the sets $\lambda S_1$ where $\lambda$ runs over
$\Lambda$. This enables us to identify the sets $\Omega$ and $S_1\times\Lambda$ so that
the adjacency matrices $X,Y,Z\in\mat_\Lambda(\ZZ S_1)$ of the relations $x$, $y$, $z$ take the
following forms:
$$
X=
\begin{pmatrix}
0 & 0 & 0 & A \\
0 & 0 & A & 0 \\
0 & uA & 0 & 0 \\
uA & 0 & 0 & 0
\end{pmatrix},
Y=
\begin{pmatrix}
0 & B & 0 & 0 \\
vB & 0 & 0 & 0 \\
0 & 0 & 0 & vB \\
0 & 0 & B & 0
\end{pmatrix},
Z=
\begin{pmatrix}
0 & 0 & wC & 0 \\
0 & 0 & 0 & wvC \\
C & 0 & 0 & 0 \\
0 & v^*C & 0 & 0
\end{pmatrix}
$$
where $A=1+a$, $B=1+b$ and $C=1+c$. Let us define the permutations $f$, $g$ and $h$ of the set
$\Omega=S_1\times\Lambda$ by means of their permutation matrices:
$$
F=
\begin{pmatrix}
0 & 0 & 0 & 1 \\
0 & 0 & 1 & 0 \\
0 & u & 0 & 0 \\
u & 0 & 0 & 0
\end{pmatrix},
G=\begin{pmatrix}
0 & 1 & 0 & 0 \\
v & 0 & 0 & 0 \\
0 & 0 & 0 & v \\
0 & 0 & 1 & 0
\end{pmatrix},
H=\begin{pmatrix}
0 & 0 & w & 0 \\
0 & 0 & 0 & wv \\
1 & 0 & 0 & 0 \\
0 & v^* & 0 & 0
\end{pmatrix}.
$$
Using the equality $uvw=1$ we obtain that $FGH=I$ (the identity matrix).
Also
$$
X=FA,\quad Y=GB,\quad Z=HC.
$$
A direct calculation shows that $F$ an $G$ commute.
Therefore $H=F^*G^*$ commute with $F$ and $G$. This implies that $F,G,H$ commute with
$X,Y,Z$. Therefore $f,g,h\in\aut(\Omega,S)$. Moreover it follows from
$F^2=uI$, $G^2=vI$, $H^2=wI$ that $\lg S_1,f,g,h\rg$ is a regular abelian group of
automorphisms of $S$. So the group $\aut(\Omega,S)$ is transitive and we are
done by Theorem~\ref{180410b}.\medskip\medskip

{\bf Acknowledgment.} The authors would like to thank Prof. M.~Hirasaka who helped
them to understand better quasi-thin schemes with at most two orthogonals and prove
the schurity of these schemes (the proof presented in our paper is different from
his proof).

The authors are also thanlful to M.~Klin for providing computational data about coherent configurations on $16$ points.

\end{document}